\documentclass[12pt]{amsart}
\usepackage{amssymb}
\usepackage{amsmath}
\usepackage{amsthm}

\newtheorem{theorem}{Theorem}[section]
\newtheorem{lemma}[theorem]{Lemma}

\theoremstyle{definition}

\newtheorem{example}[theorem]{Example}

\theoremstyle{remark}

\numberwithin{equation}{section}

\newcommand{\abs}[1]{\lvert#1\rvert}

\newcommand{\A}{{\mathcal A}}

\newcommand{\F}{{\mathcal F}}

\newcommand{\N}{\ensuremath{\mathbb{N}}}

\newcommand{\C}{\ensuremath{\mathbb{C}}}

\newcommand{\p}{\partial}

\newcommand{\Z}{\ensuremath{\mathbb{Z}}}
\begin{document}

\title[Deformation quantization of a K\"ahler-Poisson structure]{Deformation quantization of a K\"ahler-Poisson structure vanishing on a Levi nondegenerate hypersurface}

\author{Alexander V. Karabegov}
\address{Department of Mathematics, Abilene Christian University, ACU Box 28012, 215 Foster Science Building, Abilene, TX 79699-8012}
\email{alexander.karabegov@math.acu.edu}

\subjclass{Primary: 53D55}
\date{July 17, 2006}
\keywords{deformation quantization}

\begin{abstract}
We give an elementary proof of the result by Leichtnam, Tang, and Weinstein \cite{LTW} that there exists a deformation quantization with separation of variables on a complex manifold endowed with a K\"ahler-Poisson structure vanishing on a Levi nondegenerate hypersurface and nondegenerate on its complement.
\end{abstract}

\maketitle
\section{Introduction}
Let $V$ be a vector space and $\nu$ a formal parameter. Denote by $V[\nu^{-1},\nu]]$ the space of formal Laurent series of the form
\[
    v = \sum_{r \geq n} \nu^r v_r,
\]
where $n$ is possibly negative and $v_r \in V$. The elements of $V[\nu^{-1},\nu]]$ are called formal vectors.
Let $M$ be a Poisson manifold with the Poisson structure given by a Poisson bivector field $\eta$ or, equivalently, by a Poisson bracket $\{\cdot,\cdot\}$. A nondegenerate Poisson structure on $M$ is equivalent to a symplectic structure given by a symplectic form on $M$. 
Deformation quantization on $(M,\eta)$ is an associative algebra structure on $C^\infty(M)[\nu^{-1},\nu]]$ with the product $\star$ (named a star product) given by the formula
\[
    f \star g = \sum_{r \geq 0} \nu^r C_r(f,g),
\]
where $C_r$ are bidifferential operators on $M$ such that $C_0(f,g) =fg$ and $C_1(f,g) - C_1(g,f)=i\{f,g\}$.
It is assumed that the unit constant 1 is the unity in the algebra $(C^\infty(M)[\nu^{-1},\nu]],\star)$. A deformation quantization on $M$ can be localized to the open subsets of $M$. An important feature of the deformation quantizations on symplectic manifolds is the existence of local $\nu$-derivations of the form $\delta = \nu d/d\nu + A$, where $A$ is a formal differential operator (i.e., $A$ does not contain derivatives with respect to $\nu$).

Deformation quantizations were introduced in \cite{BFFLS}. The existence of deformation quantizations on the symplectic manifolds was shown by \cite{DL}, \cite{F}, and \cite{OMY}. In the general Poisson case the existence of deformation quantizations was shown by Kontsevich in \cite{K}. 

We call a complex manifold $M$ endowed with a Poisson bivector field of type (1,1) with respect to the complex structure {\it  a K\"ahler-Poisson manifold}. A nondegenerate K\"ahler-Poisson structure is equivalent to a pseudo-K\"ahler structure given by a pseudo-K\"ahler form. A deformation quantization on a K\"ahler-Poisson manifold is called deformation quantization with separation of variables (or of the Wick type) if the operators $C_r$ in the definition of the corresponding star product differentiate their first argument in antiholomorphic directions and the second argument in holomorphic ones (or vice versa). Deformation quantizations with separation of variables on an arbitrary pseudo-K\"ahler manifold were constructed in \cite{CMP1}, \cite{BW}, and \cite{RT}. It is not yet known whether there exists a deformation quantization with separation of variables on an arbitrary K\"ahler-Poisson manifold. In \cite{Third} it was shown that such quantizations cannot be obtained at least by a naive extension of the formula for the star product with separation of variables on a pseudo-K\"ahler manifold. However, it turns out that there exist 
deformation quantizations with separation of variables of a K\"ahler-Poisson structure vanishing on a Levi nondegenerate hypersurface and nondegenerate on its complement. Leichtnam, Tang, and Weinstein proved it in \cite{LTW} using para-K\"ahler Lie algebroids. Their work was motivated by the construction of Berezin-Toeplitz deformation quantization on a complex manifold with strongly pseudoconvex boundary by Engli{\v s} in \cite{E}. The goal of our paper is to give an elementary proof of this fact.

Given a star product $\star$, we denote by $L_f$ and $R_f$ the corresponding operators of left and right multiplication by an element $f$, respectively. The standard deformation quantization with separation of variables on a pseudo-K\"ahler manifold $(M,\omega)$ has the following local properties (see \cite{CMP1}). On an arbitrary contractible coordinate chart $U$ with holomorphic coordinates $\{z^k\}$
\begin{eqnarray}\label{E:ops}
  L_a = a, \ L_{\frac{1}{\nu}\frac{\p \Phi}{\p z^k}} = \frac{1}{\nu}\frac{\p \Phi}{\p z^k} + \frac{\p}{\p z^k}, \\ R_b = b, \ R_{\frac{1}{\nu}\frac{\p \Phi}{\p \bar z^l}} = \frac{1}{\nu}\frac{\p \Phi}{\p \bar z^l} + \frac{\p}{\p \bar z^l},\nonumber 
\end{eqnarray}
where 
$a$ and $b$ are holomorphic and antiholomorphic functions on $U$, respectively, and $\Phi$ is 
a potential of the pseudo-K\"ahler form $\omega$ on $U$. These properties determine the standard deformation quantization with separation of variables uniquely and globally on $M$. It was shown in \cite{LMP2} that a local $\nu$-derivation $\delta = \nu d/d\nu + A$ of the standard deformation quantization with separation of variables can be determined from the formula
\[
   \left[ e^{-\frac{1}{\nu}\Phi} \left(\nu\frac{\p}{\p \nu}\right) e^{\frac{1}{\nu}\Phi}, L_f\right] = L_{\delta(f)},
\]
where $\Phi$ is a local potential of the pseudo-K\"ahler form $\omega$.

{\bf Acknowledgments.} The author is very grateful to Professor Alan Weinstein for the inspiring discussion which stimulated this paper.

\section{A K\"ahler-Poisson structure vanishing on a Levi nondegenerate hypersurface}\label{S:som}

We begin with several calculations and reobtain some statements from \cite{LTW} in the form convenient for our exposition.
Assume that $\psi$ is a real function on a neighborhood $U \subset\C^n$. Denote  $\C^\times = \C \backslash \{0\}$ and consider the product $\tilde U = U \times \C^\times$. The holomorphic coordinates on $U$ and $\C^\times$ will be denoted by $\{z^k\}$ and $u$, respectively. Introduce a function $\rho = \psi(z,\bar z)u\bar u$ and a Hermitian matrix
\[
     \Gamma = \left( \begin{array}{cc}
\frac{\p^2 \psi}{\p z^k \p \bar z^l}u\bar u & \frac{\p \psi}{\p \bar z^l}\bar u \\
\frac{\p\psi}{\p z^k}u & \psi \end{array} \right)
\]
on $\tilde U$. We will describe the conditions on the function $\psi$ under which the matrix $\Gamma$ is nondegenerate. If this is the case, let
\[
 \Pi = \left( \begin{array}{cc}
 A^{\bar lk} & B^k \\
 C^{\bar l}& D \end{array} \right)
\]
be the inverse matrix of $\Gamma$, which is equivalent to the following conditions:
\begin{equation}\label{E:invcond1}
   \frac{\p^2 \psi}{\p z^k \p \bar z^l}u\bar u A^{\bar lm} + \frac{\p\psi}{\p z^k}u B^m = \delta_k^m, \quad \frac{\p\psi}{\p\bar z^l}\bar u A^{\bar lm} + \psi B^m = 0,
\end{equation}
\begin{equation}\label{E:invcond2}
   \frac{\p^2 \psi}{\p z^k \p \bar z^l}u\bar u C^{\bar l} + \frac{\p\psi}{\p z^k}u D = 0, \quad \frac{\p\psi}{\p\bar z^l}\bar u C^{\bar l} + \psi D = 1.
\end{equation}
The matrix $\Gamma$ can be invertible at a fixed point $\tilde x = (x,y) \in\tilde U$ (with $x\in U$ and $y \in \C^\times$) only in one of the following two cases:
\begin{enumerate}
\item $\psi(x) \neq 0$;
\item $\psi(x) = 0,\ \p\psi(x) \neq 0$, and $\bar\p\psi(x) \neq 0$.
\end{enumerate} 
Since the function $\psi$ is real, the conditions $\p\psi(x) \neq 0$ and $\bar\p\psi(x) \neq 0$ are equivalent to the condition that $x$ is not a critical point of $\psi$. We will analyse Case 1 assuming that the matrices $\Gamma$ and $\Pi$ are inverse to each other. Set
\begin{equation}\label{E:gkl}
    g_{k\bar l} = \frac{\p^2}{\p z^k \p\bar z^l} \log \abs{\psi} = \frac{1}{\psi}\left(\frac{\p^2 \psi}{\p z^k \p \bar z^l} - \frac{1}{\psi}\frac{\p\psi}{\p z^k}\frac{\p\psi}{\p\bar z^l}\right).
\end{equation}
Solving the second equation in (\ref{E:invcond1}) for $B^m$ and substituting the resulting expression to the first one, we get
\[
   \left(\frac{\p^2 \psi}{\p z^k \p \bar z^l} - \frac{1}{\psi}\frac{\p\psi}{\p z^k}\frac{\p\psi}{\p\bar z^l}\right)u\bar u A^{\bar lm} = \delta_k^m,
\]
which means, according to Eqn. (\ref{E:gkl}), that the matrix $(g_{k\bar l})$ is invertible and its inverse $(g^{\bar lk})$ is
\begin{equation}\label{E:glk}
   g^{\bar lk} = \psi u \bar u A^{\bar lk} = \rho A^{\bar lk}.
\end{equation}
Notice that $g^{\bar lk}$ does not depend on the variables $u,\bar u$ and $A^{\bar lk}$ is smooth on $\tilde U$. It easily follows from the calculations above and similar calculations applied to Eqn. (\ref{E:invcond2}) that 
in Case 1 the matrix $\Gamma$ is invertible if and only if the matrix $(g_{k\bar l})$ is invertible. 

Now consider Case 2. Let $x\in U$ be a point where $\psi(x)=0$ and which is not critical for $\psi$. Introduce the following $(n-1)$-dimensional  subspaces of $\C^n$: 
\begin{eqnarray*}
V = \left\{v = (v^1,\ldots, v^n)| \frac{\p \psi}{\p z^k}(x) v^k = 0\right\}, \\ 
W = \left\{w = (w^1,\ldots, w^n)| \frac{\p \psi}{\p \bar z^l}(x) w^{\bar l} = 0\right\}.
\end{eqnarray*}
The Levi form $Q$ is the bilinear form on $V \times W$ such that
\begin{equation}\label{E:levi}
   Q(v,w) = \frac{\p^2 \psi}{\p z^k \p \bar z^l}(x) v^k w^{\bar l}.
\end{equation}
To determine whether the matrix $\Gamma$ is invertible, consider its kernel. Assume that a nonzero vector $(v^1, \ldots v^n, a) \in \C^{n+1}$ is in the kernel of the matrix $\Gamma$, i.e.,
\begin{equation}\label{E:kernel}
    \frac{\p^2 \psi}{\p z^k \p \bar z^l}u\bar u v^k +  \frac{\p \psi}{\p \bar z^l}\bar u a = 0 \mbox{ and }  \frac{\p\psi}{\p z^k}u v^k = 0.
\end{equation}
It follows from the fisrt equation in (\ref{E:kernel}) that the vector $v=(v^1,\ldots,v^n)$ is nonzero and
$Q(v,w)=0$ for any vector $w\in W$. The second equation in (\ref{E:kernel}) means that $v=(v^1,\ldots,v^n)\in V$ and therefore the Levi form $Q$ is degenerate. Conversly, if $Q$ is degenerate, there exists a nonzero vector 
$v=(v^1,\ldots,v^n)\in V$ such that $Q(v,w)=0$ for any vector $w=(w^1,\ldots, w^n)\in W$. This means that the vector $b = (b_1, \ldots, b_n)$ such that
\[
    b_{\bar l} = \frac{\p^2 \psi}{\p z^k \p \bar z^l}(x) v^k
\]
is proportional to the vector
\[
     \left(\frac{\p\psi}{\p\bar z^1}(x), \ldots, \frac{\p\psi}{\p\bar z^l}(x)\right),
\]
which implies that the first equation in (\ref{E:kernel}) holds for some constant $a$.
Thus in Case 2 the matrix $\Gamma$ is invertible if and only if the Levi form $Q$ is nondegenerate.
Summarizing, we arrive at the following. Assume that $S \subset U$ is a Levi nondegenerate hypersurface given by a defining function $\psi$. It  means that $S$ is the zero set of the function $\psi$ which has no critical points on $S$ and the Levi form (\ref{E:levi}) is nondegenerate on $S$. The property that $S$ is Levi nondegenerate does not depend on the choice of the defining function. For any point $(x,y)\in \tilde U$ such that $x\in S$ the matrix $\Gamma$ is nondegenerate. Shrinking, if necessary, the neighborhood $U$ around $S$ we will assume from now on that $\Gamma$ is nondegenerate everywhere on $\tilde U$. In this case the function $\rho$ is a potential of the pseudo-K\"ahler form  
\[
\Omega = -i\p_{\tilde U}\bar \p_{\tilde U}\rho
\]
on $\tilde U$. Further, the function $\log \abs{\psi}$ is a potential of the pseudo-K\"ahler form
\[
    \omega = -i\p\bar\p\log \abs{\psi}
\]
and
\[
     \eta = ig^{\bar lk}\frac{\p}{\p z^k}\wedge\frac{\p}{\p\bar z^l}
\]
is the corresponding Poisson bivector field on $U\backslash S$ of type (1,1) with respect to the complex structure. Since $S$ is nowhere dense in $U$, it follows from Eqn. (\ref{E:glk}) that the Poisson bivector field $\eta$ has a smooth extension to the whole neighborhood $U$ which vanishes on $S$. Thus it determines a K\"ahler-Poisson structure on $U$ which is nondegenerate on the complement of $S$.

\section{Deformation quantization of the K\"ahler-Poisson structure on $U$}

Denote by $*$ the star product of the standard deformation quantization with separation of variables on $(\tilde U, \Omega)$ such that 
\[
    F * G = \sum_{r \geq 0} h^r  C_r(F,G).
\]
Here $h$ is a formal parameter and $F,G \in C^\infty(\tilde U)[h^{-1},h]]$. We will denote by $\tilde L_F$ and $\tilde R_G$ the left and right multiplication operators in the algebra $(C^\infty(\tilde U)[h^{-1}, h]],*)$ by the elements $F$ and $G$, respectively, so that $F * G = \tilde L_F G = \tilde R_G F$. Adapting Eqns. (\ref{E:ops}) to the star product $*$, we get
\begin{eqnarray}\label{E:tildeops}
   \tilde L_{a(z,u)} = a(z,u), \ \tilde L_{\frac{1}{h}\frac{\p\rho}{\p z^k}} = \frac{1}{h}\frac{\p\rho}{\p z^k} + \frac{\p}{\p z^k},\nonumber\\
 \tilde L_{\frac{1}{h}\frac{\p\rho}{\p u}} = \frac{1}{h}\frac{\p\rho}{\p u} + \frac{\p}{\p u},
   \tilde R_{b(\bar z,\bar u)} = b(\bar z,\bar u), \\ 
\tilde R_{\frac{1}{h}\frac{\p\rho}{\p \bar z^l}} = \frac{1}{h}\frac{\p\rho}{\p \bar z^l} + \frac{\p}{\p \bar z^l},\ \tilde R_{\frac{1}{h}\frac{\p\rho}{\p \bar u}} = \frac{1}{h}\frac{\p\rho}{\p \bar u} + \frac{\p}{\p \bar u},
\nonumber 
\end{eqnarray}
where $a(z,u)$ and $b(\bar z,\bar u)$ are a holomorphic and an antiholomorphic functions on $\tilde U$. It follows from Eqns. (\ref{E:tildeops}) that
\begin{equation}\label{E:leftrho}
   \tilde L_{\frac{1}{h}\rho} = \tilde L_{\frac{1}{h}u \frac{\p\rho}{\p u}} = \tilde L_u \tilde L_{\frac{1}{h}\frac{\p\rho}{\p u}} = \frac{1}{h} \rho + u \frac{\p}{\p u}
\end{equation}
and, similarly,
\begin{equation}\label{E:rightrho}
    \tilde R_{\frac{1}{h}\rho} = \frac{1}{h} \rho + \bar u \frac{\p}{\p \bar u}.
\end{equation}
It follows from Eqns. (\ref{E:leftrho}) and (\ref{E:leftrho}) that
\begin{equation}\label{E:inner}
   \tilde L_{\frac{1}{h}\rho} - \tilde R_{\frac{1}{h}\rho} = u \frac{\p}{\p u} - \bar u \frac{\p}{\p \bar u}.
\end{equation}
is an inner derivation of the algebra $(C^\infty(\tilde U)[h^{-1}, h]],*)$.
According to \cite{LMP2}, there exists an $h$-derivation $\delta = h\frac{d}{dh} + A$ of the star product $*$ (where $A$ is a formal differential operator on $\tilde U$) such that for  $F \in C^\infty(\tilde U)[h^{-1}, h]]$
\begin{equation}\label{E:defdelta}
   \left[e^{- \frac{1}{h}\rho}\left(h\frac{d}{dh}\right)e^{\frac{1}{h}\rho}, \tilde L_F \right] = \left[h\frac{d}{dh} - \frac{1}{h}\rho, \tilde L_F \right] = \tilde L_{\delta(F)}.
\end{equation}
Evaluating the operator in Eqn. (\ref{E:defdelta}) at the unit constant $1$ and using Eqn. (\ref{E:rightrho}), we obtain
\[
    \delta (F) = \left(h\frac{\p}{\p h} - \frac{1}{h}\rho \right) F + F * \left(\frac{1}{h}\rho \right) = \left(h\frac{\p}{\p h} + \bar u \frac{\p}{\p \bar u}\right) F.
\]
In the rest of the paper we identify the functions on $U$ with their lifts to $\tilde U$ and thus we can treat $C^\infty(U)$ as a subspace of $C^\infty(\tilde U)$. Denote by $\F(\tilde U)$ the subspace of $C^\infty(\tilde U)[h^{-1}, h]]$ consisting of the elements annihilated by the derivations (\ref{E:inner}) (i.e., commuting with $(1/h)\rho$) and $\delta$. These elements can be written as
\[
     F = \sum_{r \geq n} \left(\frac{h}{u\bar u}\right)^r f_r,
\]
where $n$ is possibly negative and $f_r \in C^\infty(U)$. The subspace $\F(\tilde U)$ is a subalgebra of the algebra $(C^\infty(\tilde U)[h^{-1}, h]],*)$ 
which does not contain the formal parameter $h$ as an element. Notice that
\[
           \frac{1}{h}\rho = \left(\frac{h}{u\bar u}\right)^{-1}\psi
\]
is a central element in this subalgebra and $C^\infty(\tilde U) \subset \F(\tilde U)$. Set $\tilde S = S \times \C^\times \subset \tilde U$. We can define an algebra 
$(\F(\tilde U \backslash \tilde S),*)$ by replacing $\tilde U$ with $\tilde U \backslash \tilde S$ in the construction of $\F(\tilde U)$. Since $\rho$ does not vanish on $\tilde U \backslash \tilde S$, the element $(1/h)\rho$ has an inverse in the algebra $(\F(\tilde U \backslash \tilde S),*)$ which is also central. Denote it by $\kappa$. Thus $\kappa = h/\rho \pmod{h^2}$ and 
\[
    \tilde L_\kappa = \left(\tilde L_{\frac{1}{h}\rho}\right)^{-1}.
\]
For any integer $n$, denote by $\kappa^{*n}$ the $n$th power of the element $\kappa$ with respect to the star product $*$. In particular, $\kappa^{*(-1)} = (1/h)\rho$. Since the operator $\tilde L_{(1/h)\rho}$ commutes with the pointwise multiplication by the functions from $C^\infty(U\backslash S)$, so does $\tilde L_\kappa$, and therefore for any $f \in C^\infty(U \backslash S)$ and $n \in \Z$
\[
   f * \left (\kappa^{*n}\right) = \left (\kappa^{*n}\right) * f = \left(\tilde L_\kappa\right)^n f = f\cdot \left(\tilde L_\kappa\right)^n 1 =  \left (\kappa^{*n}\right) \cdot f.
\]
Let $\nu$ be a different formal parameter. Define a mapping 
\[
\tau: C^\infty(U \backslash S)[\nu^{-1},\nu]] \to \F(\tilde U \backslash \tilde S)
\]
via an $h$-adically convergent series:
\[
   \tau \left(\sum_{r\geq n}\nu^r f_r(z,\bar z)\right) = \sum_{r\geq n} \kappa^{*r} *  f_r(z,\bar z) = \sum_{r\geq n} \kappa^{*r}\cdot  f_r(z,\bar z).
\]
Notice that 
\begin{eqnarray}\label{E:cinftymodule}
  \tau(f) = f, \ \tau (F \cdot f) = \tau(F) \cdot f,\ \tau(\nu) = \kappa, \\
 \mbox{ and } \tau (\nu F) = \kappa * \tau(F) = \tau(F) * \kappa   \nonumber  
\end{eqnarray}
for any $f \in C^\infty(U \backslash S)$ and $F \in \F(\tilde U \backslash \tilde S)$. Using the fact that $\kappa^{*n} = (h/\psi)^n \pmod{h^{n+1}}$, it is easy to show that $\tau$ is a bijection of the space $C^\infty(U \backslash S)[\nu^{-1},\nu]]$ onto $\F(\tilde U \backslash \tilde S)$. Denote by $\star$ the pullback of the star-product $*$ to $C^\infty(U \backslash S)[\nu^{-1},\nu]]$ via the mapping $\tau$ and by $L_f$ and $R_f$ the operators of left and right multiplication by $f \in C^\infty(U \backslash S)[\nu^{-1},\nu]]$  with respect to the product $\star$. One can show with the use of Eqns. (\ref{E:cinftymodule}) that $\star$ is actually a star product. Assume that $a$ and $b$ are a holomorphic and an antiholomorphic functions on $U\backslash S$, respectively. Then, using Eqns. (\ref{E:cinftymodule}), we get for $f \in C^\infty(U \backslash S)[\nu^{-1},\nu]]$ that
\[
   \tau(a \star f) = \tau(a) * \tau(f) = a \cdot \tau(f) = \tau(a\cdot f),
\]
which means that $a \star f = af$. Similarly, one can show that $f \star b = f\cdot b$. Thus $\star$ is a star product of a deformation quantization with separation of variables on $U\backslash S$.

Consider the operator
\[
    A_k = \frac{1}{\nu} \frac{\p}{\p z^k}\log \abs{\psi}
 +  \frac{\p}{\p z^k} = \frac{1}{\nu \psi} \frac{\p \psi}{\p z^k}+ \frac{\p}{\p z^k}
\]
on $C^\infty(U\backslash S)[\nu^{-1},\nu]]$. We will need the following technical statement.
\begin{lemma}\label{L:technstat}
Given a function $f \in C^\infty(U\backslash S)$, the formula
\[
    \left(\frac{1}{h}\frac{\p \rho}{\p z^k}\right) * f = \tau(A_kf).   
\]
holds.
\end{lemma}
{\it Proof.} Using Eqns (\ref{E:tildeops}), we get
\begin{eqnarray*}
   \left(\frac{1}{h}\frac{\p \rho}{\p z^k}\right) * f = \left(\frac{1}{h}\frac{\p \rho}{\p z^k}\right)f + \frac{\p f}{\p z^k} = \\
\left(\frac{1}{h}\rho\frac{\p}{\p z^k}\log \abs{\psi}
  \right)f +  \frac{\p f}{\p z^k} =
\frac{1}{h}\rho * \left( \left(\frac{\p}{\p z^k}\log \abs{\psi}
 \right)f\right) +\\
 \frac{\p f}{\p z^k} = 
\tau\left( \frac{1}{\nu}\left(\frac{\p}{\p z^k}\log \abs{\psi}
 \right)f + \frac{\p f}{\p z^k}\right) = \tau(A_kf),
\end{eqnarray*}
which concludes the proof.

\begin{lemma}\label{L:leftmult}
 The operator of left multiplication by $\frac{1}{\nu} \frac{\p}{\p z^k}\log \abs{\psi}
 $ with respect to the product $\star$ coincides with the operator $A_k$,
\[
    L_{\frac{1}{\nu} \frac{\p}{\p z^k}\log \abs{\psi}
 } = A_k.
\]
\end{lemma}
{\it Proof.} Given an element $f = \sum_{r \geq n} \nu^f f_r \in C^\infty(U\backslash S)[\nu^{-1},\nu]]$,
we get from Eqns. (\ref{E:tildeops}), Lemma \ref{L:technstat}, and the fact that $\kappa$ is central, that
\begin{eqnarray*}
  \tau \left(\frac{1}{\nu} \frac{\p}{\p z^k}\log \abs{\psi}
 \star f\right) = \left(\frac{1}{h}\rho \frac{\p}{\p z^k}\log \abs{\psi}
 \right) * \tau (f) =\\
 \left(\frac{1}{h}\frac{\p \rho}{\p z^k}\right) * \tau (f) = 
\sum_{r \geq n} \kappa^{*r} * \left(\frac{1}{h}\frac{\p \rho}{\p z^k}\right) * f_r =\\
  \sum_{r \geq n} \kappa^{*r} * \tau(A_kf_r) =
\sum_{r \geq n} \tau(\nu^r A_kf_r) =\tau (A_kf),
\end{eqnarray*}
which proves the Lemma. 

Since $\log \abs{\psi}$ is a potential of the pseudo-K\"ahler form $\omega$ on $U\backslash S$, Lemma \ref{L:leftmult} immediately implies the following
\begin{theorem}\label{T:ident}
The product $\star$ is the star product of the standard deformation quantization with separation of variables on $(U \backslash S, \omega)$.
\end{theorem}

We want to show that the star product $\star$ can be extended to the whole neighborhood $U$. This requires some preparations.
For $r \in \N$ denote by $N_r(\nu)$ the following formal number:
\[
    N_r(\nu) = \prod_{s = 1}^r \frac{\nu}{1 + \nu s} = \nu^r + \ldots,
\]
where the division by $1 + \nu s$ is in the field of formal numbers.
Set $N_0(\nu) = 1$.
\begin{lemma}\label{L:invtau}
  Given a function $f \in C^\infty(U\backslash S)$, the following formula holds for any integer $r \geq 0$:
\[
   \tau^{-1}\left(\left(\frac{h}{u\bar u}\right)^rf\right) = N_r(\nu) \psi^r f.
\]
\end{lemma}
{\it Proof.} Using Eqns. (\ref{E:tildeops}), we have for any $s \in \N$:
\begin{eqnarray*}
  \frac{1}{\nu} \tau^{-1}\left(\left(\frac{h}{u\bar u}\right)^sf\right) = \tau^{-1}\left(\frac{1}{h}\rho\right) \star \tau^{-1}\left(\left(\frac{h}{u\bar u}\right)^sf\right) =\\
\tau^{-1}\left(\frac{1}{h}\rho *\left(\frac{h}{u\bar u}\right)^sf\right)=
\tau^{-1}\left(\frac{1}{h}\rho \left(\frac{h}{u\bar u}\right)^sf + u\frac{\p}{\p u}\left(\left(\frac{h}{u\bar u}\right)^sf\right)\right)= \\
\tau^{-1}\left(\left(\frac{h}{u\bar u}\right)^{s-1}\psi f\right) - s \tau^{-1}\left(\left(\frac{h}{u\bar u}\right)^sf\right),
\end{eqnarray*}
whence
\[
   \frac{1+ \nu s}{\nu} \tau^{-1}\left(\left(\frac{h}{u\bar u}\right)^sf\right) = \tau^{-1}\left(\left(\frac{h}{u\bar u}\right)^{s-1}\psi f\right).
\]
Now the statement of the Lemma follows by induction.

Denote by $\F_{\geq 0}(\tilde U)$ the subspace of $\F(\tilde U)$
consisting of the elements
\[
   f=\sum_{r\geq 0}\left(\frac{h}{u \bar u}\right)^r f_r,
\]
where $f_r \in C^\infty(U)$. It is a subalgebra of $(\F(\tilde U),*)$. Define a mapping $\sigma : \F_{\geq 0}(\tilde U) \to C^\infty(U)[[\nu]]$ by the formula
\[
   \sigma\left(\sum_{r\geq 0}\left(\frac{h}{u \bar u}\right)^r f_r\right) = \sum_{r\geq 0} N_r(\nu)\psi^r f_r.
\]
Denote $\A(U) = \sigma(\F_{\geq 0}(\tilde U)) \subset C^\infty(U)[[\nu]]$. Since the mapping $\sigma$ is injective and local, one can push forward the star product $*$ via $\sigma$ to a product on $\A(U)$ localizable to the open subsets of $U$ (where it will be defined on the restrictions of the elements from $\A(U)$). It follows from Lemma \ref{L:invtau} that the restriction of that product to $U\backslash S$ agrees with the product $\star$. Thus it will be denoted by $\star$ as well. Notice that since $C^\infty(U) \subset \F_{\geq 0}(\tilde U)$, the product $f * g$ of elements $f,g \in C^\infty(U)$ belongs to $\F_{\geq 0}(\tilde U)$ and thus the term $h^r  C_r(f,g)$ of the formal series representing $f * g$ can be rewritten as
\[
    h^r  C_r(f,g) = \left(\frac{h}{u\bar u}\right)^r D_r(f,g) 
\]
for some bidifferential operator $D_r$ on $U$. Since $C^\infty(U) \subset \A(U)$ and $\sigma(f) = f$ for any $f \in C^\infty(U)$, we see that for $f,g \in C^\infty(U)$
\begin{eqnarray*}
   f \star g = \sigma(f * g) = \sigma\left(\sum_{r\geq 0} h^r  C_r(f,g)\right) =\\
 \sigma\left(\sum_{r\geq 0} \left(\frac{h}{u\bar u}\right)^r D_r(f,g)\right) =
\sum_{r\geq 0} N_r(\nu)\psi^r D_r(f,g).
\end{eqnarray*}
This formula gives a smooth extension of the star product $\star$ from $U \backslash S$ to $U$ and shows that $f \star g - fg = 0$ on $S$. The results of this paper can be globalized as follows. Assume that $(M,\eta)$ is a K\"ahler-Poisson  manifold such that
\begin{itemize}
\item the K\"ahler-Poisson bivector field $\eta$ is nondegenerate on the complement of a Levi-nondegenerate hypersurface $S \subset M$, i.e., $\eta$ determines a pseudo-K\"ahler form $\omega$ on $M\backslash S$; and
 \item for any point $x\in S$ there exists a local defining function $\psi$ of the hypersurface $S$ in a neighborhood of $x$ such that $\log \abs{\psi}$ is a potential of the form $\omega$ on the complement of $S$.
\end{itemize}
The results obtained in this paper imply the following theorem 
\begin{theorem}
  The standard deformation quantization with separation of variables on the pseudo-K\"ahler manifold 
$(M\backslash S,\omega)$ extends to a deformation quantization with separation of variables on the K\"ahler-Poisson manifold $(M,\eta)$.
\end{theorem}

\begin{example}
Consider a defining function 
\[
   \psi = \sum_{k = 1}^n \abs{z^k}^2 - 1
\]
of the unit sphere $S \subset \C^n$. The unit sphere $S$ is a Levi nondegenerate hypersurface in $\C^n$.
On the complement of $S$ the potential $\log\abs{\psi}$ determines the pseudo-K\"ahler metric
\[
      g_{kl} = \frac{1}{\psi}\left(\delta_{kl} - \frac{1}{\psi}\bar z^k z^l\right).
\] 
Its inverse
\[
     g^{lk} = \psi \left(\delta^{kl} - \bar z^l z^k\right)
\]
is a K\"ahler-Poisson tensor which gives a global  K\"ahler-Poisson structure on $\C^n$. It extends via the inclusion $\C^n \subset \C P^n$ given by the formula
\[
    (z^1, \ldots, z^n) \mapsto (z^1: \ldots :z^n:1)
\]
to a K\"ahler-Poisson structure on the complex projective space $\C P^n$ invariant with respect to the projective action of the group $SU(n,1)$. There exists a global $SU(n,1)$-invariant star product on $\C P^n$ which coincides with the star product of the standard deformation quantization with separation of variables on the complement of the hypersurface $S \subset \C^n \subset \C P^n$. This invariant star product can be constructed also by the methods developed in \cite{AL}. 
\end{example}

\bibliographystyle{amsalpha}

\end{document}